\documentclass{llncs}

\usepackage{graphicx}      
\usepackage{subfig}

\usepackage{amsmath} 

\usepackage{amssymb} 
\usepackage{commath}

\newtheorem{thm}{Theorem}

\newcommand{\G}{G}
\newcommand{\g}{\mathfrak{g}}
\newcommand{\h}{\mathfrak{h}}
\newcommand{\m}{\mathfrak{m}}
\newcommand{\R}{\mathbb{R}}
\newcommand{\cP}{\mathcal{P}}
\DeclareMathOperator{\Ad}{Ad}

\begin{document}

\title{Shape analysis on Lie groups and homogeneous spaces}

\author{Elena Celledoni\inst{1} \and S\o{}lve Eidnes\inst{1}, Markus Eslitzbichler\inst{2} \and Alexander Schmeding\inst{1}} 
\authorrunning{Celledoni, Eidnes, Eslitzbichler, Schmeding} 
\tocauthor{E. Celledoni, S. Eidnes, M. Eslitzbichler and A. Schmeding} 
 \institute{NTNU Trondheim, Trondheim, Norway,\\ \email{elena.celledoni@ntnu.no},\\ \email{solve.eidnes@ntnu.no},\\ 
 \email{m.eslitzbichler@gmail.com}, \\ \email{schmeding@tu-berlin.de} }

\maketitle              

\begin{abstract}

In this paper we are concerned with the 
approach to shape analysis based on the so called Square Root Velocity Transform (SRVT). 
We propose a generalisation of the SRVT from Euclidean spaces to shape spaces of curves on Lie groups and on homogeneous manifolds. The main idea behind our approach is to exploit the geometry of the natural Lie group actions on these spaces.

 \end{abstract}  

Shape analysis methods have significantly increased in popularity in the last decade. Advances in this field have been made both in the theoretical foundations and in the extension of the methods to new areas of application. Originally developed for planar curves, the techniques of shape analysis have been successfully extended to higher dimensional curves, surfaces, activities, character motions and a number of different types of digitalized objects. 

In the present paper, shapes are unparametrized curves, evolving on a vector space, on a Lie group, or on a manifold. Shape spaces and spaces of curves are infinite-dimensional Riemannian manifolds, whose Riemannian metrics are the crucial tool to compare and analyse shapes. 

We are concerned with one particular approach to shape analysis, which is based on the Square Root Velocity Transform (SRVT) \cite{srivastava11sao}. On vector spaces, the SRVT maps parametrized curves (i.e.\ smooth immersions) to appropriately scaled tangent vector fields along them via 
\begin{equation}\label{EUC: SRVT}
\mathcal{R} \colon \mathrm{Imm} ([0,1], \R^d) \rightarrow C^\infty ([0,1], \R^d \setminus \{0\}), \quad c \mapsto \frac{\dot{c}}{\sqrt{\norm{\dot{c}}}}.
\end{equation}
The transformed curves are then compared computing geodesics in the $L^2$ metric, and the scaling induces reparametrization invariance of the pullback metric.  
Note that it is quite natural to consider an $L^2$ metric directly on the original parametrized curves. Constructing the $L^2$ metric with respect to integration by arc-length, one obtains a reparametrisation invariant metric. However, this metric is unsuitable for our purpose as it leads to vanishing geodesic distance on the quotient shape space \cite{michor05vgd} and consequently also on the space of parametrised curves \cite{MR2891297}. This infinite-dimensional phenomenon prompted the investigation of alternative, higher order Sobolev type metrics \cite{michor06rgo}, which however can be computationally demanding. 
Since 
%
it allows geodesic computations via the $L^2$ metric on the transformed curves, the  SRVT technique is computationally attractive. 
It is also possible to prove that this algorithmic approach corresponds, at least locally, to a particular Sobolev type metric, see \cite{bauer14cri,celledoni15sao}.

We propose a generalisation of the SRVT to construct well-behaved Riemannian metrics on shape spaces with values in Lie groups and homogeneous manifolds.
Our methodology is alternative to what was earlier proposed in \cite{su14sao,brigant} and the main idea is, following \cite{celledoni15sao}, to take advantage of the Lie group acting transitively on the homogeneous manifold. Since we want to compare curves, the main tool here is an SRVT which transports the manifold valued curves into the Lie algebra or a  subspace of the Lie algebra. 

\section{SRVT for Lie group valued shape spaces}
\begin{figure}[htbp] 
	\begin{center}
		\includegraphics[width=.9\textwidth]{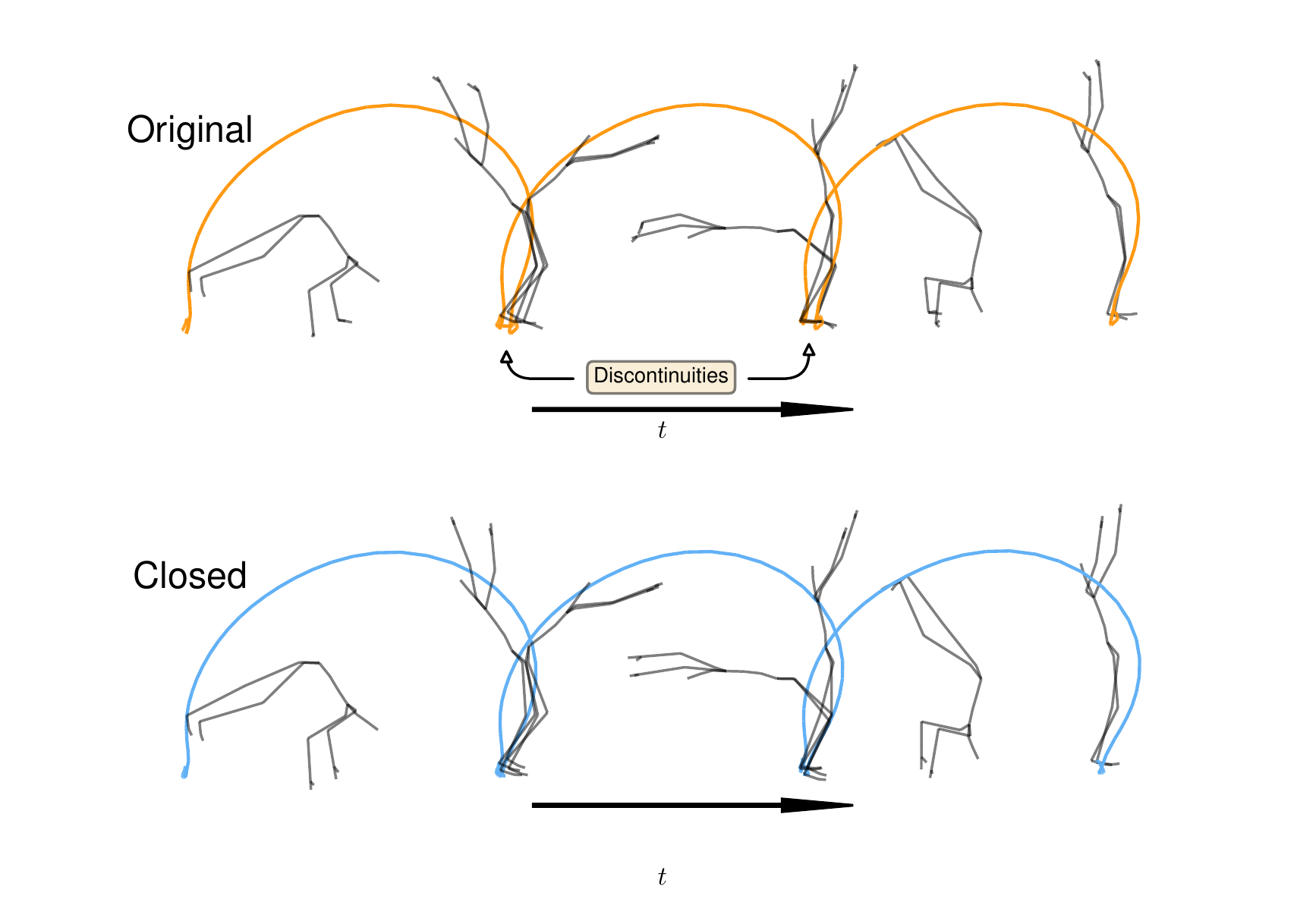}
	\end{center}
	\caption{Application of closing algorithm to a cartwheel animation.  
		Note the large difference between start and end poses, on the right and the left respectively.
		The motion is repeated once and suffers from a strong jerk when it repeats, especially in the left hand.
		In the second row, the curve closing method has been used to alleviate this discontinuity.
	}\label{fig1}
\end{figure}

In the Lie group case, the obvious choice for this tangent space is of course the Lie algebra $\g$ of the Lie group $\G$. The idea is to use the derivative $T_e R_g$ of the right translation for the transport and measure with respect to a right-invariant Riemannian metric.\footnote{Equivalently one could instead use left translations and a left-invariant metric here.} Instead of the ordinary derivative, one thus works with the right-logarithmic derivative $\delta^r (c) (t) = T_e R_{c(t)^{-1}} (\dot{c}(t))$ (here $e$ is the identity element of $\G$) and defines an SRVT for Lie group valued curves as (see \cite{celledoni15sao}):
\begin{equation}\label{LG: SRVT}
\mathcal{R} \colon \mathrm{Imm} ([0,1], \G) \rightarrow C^\infty ([0,1] , \g \setminus \{0\}), \quad c \mapsto \frac{\delta^r (c)}{\sqrt{\norm{\dot{c}}}}.
\end{equation} 
We will use the short notetion $I=[0,1]$ in what follows. Using tools from Lie theory, we are then able to describe the resulting pullback metric on the space $\cP_*$ of immersions $ c \colon [0,1] \rightarrow G$ which satisfy $c(0) = e$: 

\begin{thm}[The Elastic metric on Lie group valued shape spaces {\cite{celledoni15sao}}]
	Let $c \in \cP_*$ and consider $v,w \in T_c \cP_*$.
	The pullback of the $L^2$-metric on $C^\infty(I, \g \setminus \{0\})$ under the SRVT \eqref{LG: SRVT} to $\cP_*$ is given by the first order Sobolev metric: 
	\begin{equation}\label{Eq:ElasticMetric} \begin{aligned}
	G_c (v,w) = \int_I \frac{1}{4} &\left\langle D_s v, u_c\right\rangle \left\langle D_s w, u_c\right\rangle  \\ &+ \left\langle D_s v-u_c \left\langle D_s v,u_c \right\rangle , D_s w-u_c\left\langle D_s w,u_c \right\rangle \right\rangle  \dif s,
	\end{aligned}
	\end{equation}
	where $D_s v := T_c \delta^r (v)/ \norm{ \dot{c} }$, $u_c := \delta^r(c)/\norm{\delta^r (c)}$ is the unit tangent vector of $\delta^r (c)$ and $\dif s = \norm{\dot{c}(t)} \dif t$.
\end{thm}
 The geodesic distance of this metric descends to a nonvanishing metric on the space of unparametrized curves. In particular, this distance is easy to compute as one can prove \cite[Theorem 3.16]{celledoni15sao} that 
 \begin{thm}
  If $\text{\upshape dim } \g > 2$, then the geodesic distance of $C^\infty (I,\g\setminus \{0\})$ is globally given by the $L^2$-distance.  
  In particular, in this case the geodesic distance of the pullback metric \eqref{Eq:ElasticMetric} on $\cP_*$ is given by 
  $$  d_{\cP_*}(c_0, c_1) := \sqrt{\int_I \| \mathcal{R} (c_0) (t) - \mathcal{R} (c_1)(t) \|^2 \dif t}.$$
 \end{thm}	
These tools give rise to algorithms which can be used in, among other things, tasks related to computer animation and blending of curves, as shown in \cite{celledoni15sao}.
The blending $c(t,s)$ of two curves $c_0(t)$ and $c_1(t)$, $t\in I$, amounts simply to a convex linear convex combination of their SRV transforms:
$$c(t,s)=\mathcal{R}^{-1}\left(s\, \mathcal{R}(c_0(t))+ (1-s) \mathcal{R}(c_1(t))\right),\qquad s \in [0,1].$$
Using the transformation of the curves to the Lie algebra by the SRVT, we also propose a curve closing algorithm allowing one to remove discontinuities from motion capturing data while preserving the general structure of the movement. (See Figure \ref{fig1}.)

\section{The structure of the SRVT}
Analysing the constructions for the square root velocity transform, e.g.\ \eqref{EUC: SRVT} and \eqref{LG: SRVT} or the generalisations proposed in the literature, every SRVT is composed of three distinct building blocks. While two of these blocks can not be changed, there are many choices for the second one (transport) in constructing an SRVT:
\begin{itemize}
	\item \textbf{Differentiation}: The basic building block of every SRVT, taking a curve to its derivative.
	\item \textbf{Transport}: Bringing a curve into a common space of reference. In general there are many choices for this transport\footnote{In the literature, e.g.\ \cite{su14sao}, a common choice is parallel transport with respect to the Riemannian structure.} (in our approach we use the Lie group action to transport data into the Lie algebra of the acting group).
	\item \textbf{Scaling}: The second basic building block, assures reparametrization invariance of the metrics obtained.
\end{itemize}
In constructing the SRVT, we advocate the use of Lie group actions for the transport. This action allows us to transport derivatives of curves to our choice of base point and to lift this information to a curve in the Lie algebra.

Other common choices for the transport usually arise from parallel transport (cf.\ e.g.\ \cite{su14sao,brigant}). 
The advantage of using the Lie group action is that we obtain a global transport, i.e.\ we do not need to restrict to certain open submanifolds to make sense of the (parallel) transport.\footnote{The problem in these approaches arises from choosing curves along which the parallel transport is conducted. Typically, one wants to transport along geodesics to a reference point and this is only well-defined outside of the cut locus (also cf.\ \cite{1612.02604v1}).} Last but not least, right translation is in general computationally more efficient than computing parallel transport using the original Riemannian metric on the manifold.

\section{SRVT on homogeneous spaces} 
 
Our approach \cite{1704.01471v1} for shape analysis on a homogeneous manifold $\mathcal{M} = \G / H$ exploits again the geometry induced by the canonical group action $\Lambda \colon G \times \mathcal{M} \rightarrow \mathcal{M}$. We fix a Riemannian metric on $\G$ which is right $H$-invariant, i.e. the maps $R_h$ for $h\in H$ are Riemannian isometries.
The SRVT is obtained using a right inverse of the composition of the Lie group action with the evolution operator (i.e.\ the inverse of the right-logarithmic derivative) of the Lie group. 
If the homogeneous manifold is reductive,\footnote{Recall that a homogeneous space $\G/H$ is reductive if the Lie subalgebra $\h$ of $H \subseteq \G$ admits a reductive complement, i.e.\ $\g = \h \oplus \m$, where $\m$ is a subvector space invariant under the adjoint action of $H$.} there is an explicit way to construct this right inverse. Identifying the tangent space at $[e]$, the equivalence class of the identity, via $\omega_e \colon T_{[e]} \mathcal{M} \rightarrow \m \subseteq \g$ with the reductive complement.
Then we define the map $\omega ([g]) = \Ad (g). \omega_e (T\Lambda (g^{-1} ,\cdot) [g])$ (which is well-defined by reductivity) and obtain a square root velocity transform for reductive homogeneous spaces as 
\begin{equation} \label{SRVT hom}
                  \mathcal{R} \colon \mathrm{Imm} ([0,1], \mathcal{M}) \rightarrow C^\infty ([0,1], \g \setminus \{0\}) ,\quad  c\mapsto \frac{ \omega \circ \dot{c}}{\sqrt{\norm{\omega \circ \dot{c}}}}
\end{equation}
Conceptually this SRVT is somewhat different from the one for Lie groups, as it does not establish a bijection between the manifolds of smooth mappings.
However, one can still use \eqref{SRVT hom} to construct a pullback metric on the manifold of curves to the homogeneous space by pulling back the $L^2$ inner product of curves on the Lie algebra through the SRVT. 
Different choices of Lie group actions will give rise to different Riemannian metrics (with different properties). 

\section{Numerical experiments}

We present some results about the realisation of this metric through the SRVT framework in the case of reductive homogeneous spaces.
Further, our results are illustrated in a concrete example. 
We compare the new methods for curves into the sphere $\mathrm{SO} (3) / \mathrm{SO} (2)$ with results derived from the Lie group case.

In the following, we use the Rodrigues' formula for the Lie group exponential $\exp \colon \mathfrak{so}(3) \rightarrow \mathrm{SO} (3)$,
\begin{equation*}
\exp(\hat{x}) = I + \frac{\sin{(\alpha)}}{\alpha}\hat{x} + \frac{1- \cos{(\alpha)}}{\alpha^2}\hat{x}^2, \qquad \alpha = \lVert x \rVert_2
\end{equation*}
and the corresponding formula for the logarithm $\log \colon \mathrm{SO}(3) \rightarrow \mathfrak{so}(3)$,
\begin{equation*}
\log(X) = \frac{\sin^{-1}(\lVert y \rVert)}{\lVert y \rVert}\hat{y}, \quad X \neq I, \quad X \text{ close to } I,
\end{equation*}
are used, where $\hat{y} = \frac{1}{2}(X-X^\text{T})$, and the relationship between $x$ and $\hat{x}$ is given by the isomorphism between $\mathbb{R}^3$ and $\mathfrak{so}(3)$ known as the hat map
\begin{equation*}
x = 
 \begin{pmatrix}
  x_1 \\
  x_2 \\
  x_3
 \end{pmatrix}
\mapsto
\hat{x} = 
 \begin{pmatrix}
  0 & -x_3 & x_2 \\
  x_3 & 0 & -x_1 \\
  -x_2 & x_1 & 0
 \end{pmatrix}.
\end{equation*}

\subsection{Lie group case} Consider a continuous curve $z(t), t \in [t_0,t_N]$, in $\mathrm{SO}(3)$. We approximate it by $\bar{z}(t)$, interpolating between $N+1$ values $\bar{z}_i = z(t_i)$, with $t_0 < t_1 < ... < t_N$, as:
\begin{equation}
\bar{z}(t) := \sum_{i=0}^{N-1} \chi_{[t_i,t_{i+1})}(t) \exp{\left(\frac{t-t_i}{t_{i+1}-t_i} \log\left(\bar{z}_{i+1} \bar{z}_{i}^\text{T} \right)\right)}\bar{z}_i,
\label{eq:zapprox}
\end{equation}
where $\chi$ is the characteristic function.

The SRVT (\ref{LG: SRVT}) of $\bar{z}(t)$ is a piecewise constant function $\bar{p}(t)$ in $\mathfrak{so}(3)$ with values $\bar{p}_i = \bar{p}(t_i), \,i = 0,...,N-1$, found by
\begin{align*}
\bar{p}_i = \frac{\eta_i}{\sqrt{\lVert \eta_i \rVert}}, \qquad \eta_i = \frac{\log(\bar{z}_{i+1}\bar{z}_i^\text{T})}{t_{i+1}-t_i}.
\end{align*}
The inverse $\mathcal{R}^{-1} : \mathfrak{so}(3) \rightarrow \mathrm{SO}(3)$ is then given by (\ref{eq:zapprox}), with the discrete points
\begin{equation*}
\bar{z}_{i+1} = \exp{\left(\lVert\bar{p}_i\rVert\bar{p}_i\right)}\bar{z}_i, \quad i = 1,...,N-1, \quad \bar{z}_0 = z(t_0).
\end{equation*}

\subsection{Homogeneous manifold case} As an example of the homogeneous space case, consider the curve $c(t)$ on the sphere SO(3)/SO(2) (i.e. S$^2$), which we approximate by $\bar{c}(t)$, interpolating between the $N+1$ values $\bar{c}_i =c(t_i)$:
\begin{equation}
\bar{c}(t) := \sum_{i=0}^{N-1} \chi_{[t_i,t_{i+1})}(t) \exp{\left(\frac{t-t_i}{t_{i+1}-t_i}\left(v_i \bar{c}_i^\text{T}- \bar{c}_i v_i^\text{T} \right)\right)}\bar{c}_i,
\label{eq:capprox}
\end{equation}
where $v_i$ are approximations to $\left.\frac{d}{d t}\right|_{t=t_i}c(t)$ found by solving the equations
\begin{align}
\label{eqforvi}
\bar{c}_{i+1} = \exp{\left(v_i \bar{c}_i^\text{T}- \bar{c}_i v_i^\text{T}\right)}\bar{c}_i, \\
\text{constrained by } \quad v_i^\text{T}\bar{c}_i = 0.
\label{eq:constrains}
\end{align}
Observing that if $\kappa = \bar{c}_i \times v_i$, then $\hat{\kappa} = v_i \bar{c}_i^\text{T}- \bar{c}_i v_i^\text{T}$, and assuming that the sphere has radius $1$, we have by (\ref{eq:constrains}) that $\lVert \bar{c}_i \times v_i \rVert_2 = \lVert \bar{c}_i \rVert_2 \lVert v_i \rVert_2 = \lVert v_i \rVert_2$. By \eqref{eqforvi} we get 
$$
\bar{c}_{i+1} = 
\frac{\sin{(\lVert v_i \rVert_2)}}{\lVert v_i \rVert_2}v_i + \cos{\left(\lVert v_i \rVert_2\right)} \bar{c}_i.
$$
Calculations give $\bar{c}_i^\text{T}\bar{c}_{i+1} = 1- \cos{\left(\lVert v_i \rVert_2\right)}$ and $\lVert v_i \rVert_2 = \arccos{\left(\bar{c}_i^\text{T}\bar{c}_{i+1}\right)}$, leading to
$v_i = \left(\bar{c}_{i+1} - \bar{c}_i^\text{T}\bar{c}_{i+1}\bar{c}_i\right) \frac{\arccos{\left(\bar{c}_i^\text{T}\bar{c}_{i+1}\right)}}{\sqrt{1-\left(\bar{c}_i^\text{T}\bar{c}_{i+1}\right)^2}},$ which we insert into (\ref{eq:capprox}) to get
\begin{equation}
\bar{c}(t) = \sum_{i=0}^{N-1} \chi_{[t_i,t_{i+1})}(t) \exp{\left(\frac{t-t_i}{t_{i+1}-t_i} \frac{\arccos{\left(\bar{c}_i^\text{T}\bar{c}_{i+1}\right)}}{\sqrt{1-\left(\bar{c}_i^\text{T}\bar{c}_{i+1}\right)^2}} \left(\bar{c}_{i+1} \bar{c}_i^\text{T} - \bar{c}_i \bar{c}_{i+1}^\text{T} \right)\right)}\bar{c}_i.
\label{eq:capprox2}
\end{equation}
The SRVT (\ref{SRVT hom}) of $\bar{c}(t)$ is a piecewise constant function $\bar{q}(t)$ in $\mathfrak{so}(3)$, taking values $\bar{q}_i = \bar{q}(t_i), \,i = 0,...,N-1$, where
\begin{align*}
\bar{q}_i &= \mathcal{R}(\bar{c}_i) = \frac{a_{\bar{c}_i}(v_i)}{\lVert a_{\bar{c}_i}(v_i) \rVert^\frac{1}{2}} = \frac{v_i \bar{c}_i^\text{T}- \bar{c}_i v_i^\text{T}}{\lVert v_i \bar{c}_i^\text{T}- \bar{c}_i v_i^\text{T} \rVert^\frac{1}{2}}\\
&= \frac{\arccos^\frac{1}{2}{\left(\bar{c}_i^\text{T}\bar{c}_{i+1}\right)}}{\left(1-\left(\bar{c}_i^\text{T}\bar{c}_{i+1}\right)^2\right)^\frac{1}{4}\lVert \bar{c}_{i+1} \bar{c}_i^\text{T} - \bar{c}_i \bar{c}_{i+1}^\text{T}  \rVert^\frac{1}{2}} \left(\bar{c}_{i+1} \bar{c}_i^\text{T} - \bar{c}_i \bar{c}_{i+1}^\text{T} \right)
\end{align*}
The inverse of this SRVT is given by (\ref{eq:capprox2}), with the discrete points found as in the Lie group case by $\bar{c}_{i+1} = \exp{\left(\lVert\bar{q}_i\rVert\bar{q}_i\right)}\bar{c}_i$ and $\bar{c}_0 = c(t_0)$.

As an alternative, we define the reductive SRVT \cite{1704.01471v1} by $$\mathcal{R}_{\mathfrak{m}}(\bar{c}_i) := \mathcal{R}([U,U^\perp]^\text{T}_i \bar{c}_i),$$
where $[U,U^\perp]_{i+1} = \exp{(a_{\bar{c}_i}(v_i))}[U,U^\perp]_{i}$ for $i=0,...,N-1$, and $[U,U^\perp]_0$ can be found e.g. by $QR$-factorization of $c(t_0)$.

\begin{center}
\begin{figure}[htbp]
\subfloat[From left to right: Two curves on the sphere, their original parametrizations, the reparametrization minimizing the distance in $\text{SO}(3)$ and the reparametrization minimizing the distance in $S^\text{2}$, using the reductive SRVT.]{
\includegraphics[width=0.24\textwidth]{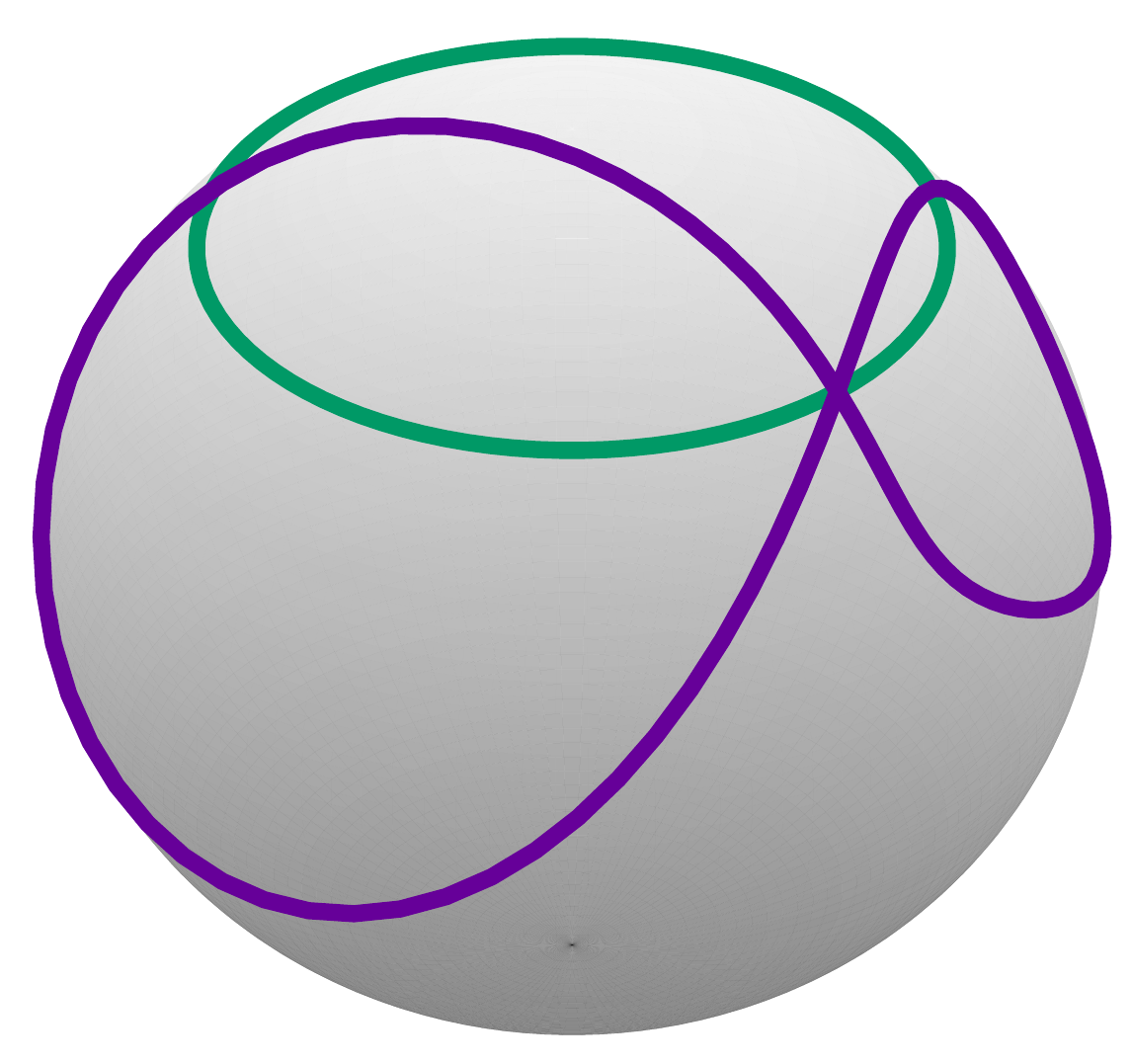}
\includegraphics[width=0.24\textwidth]{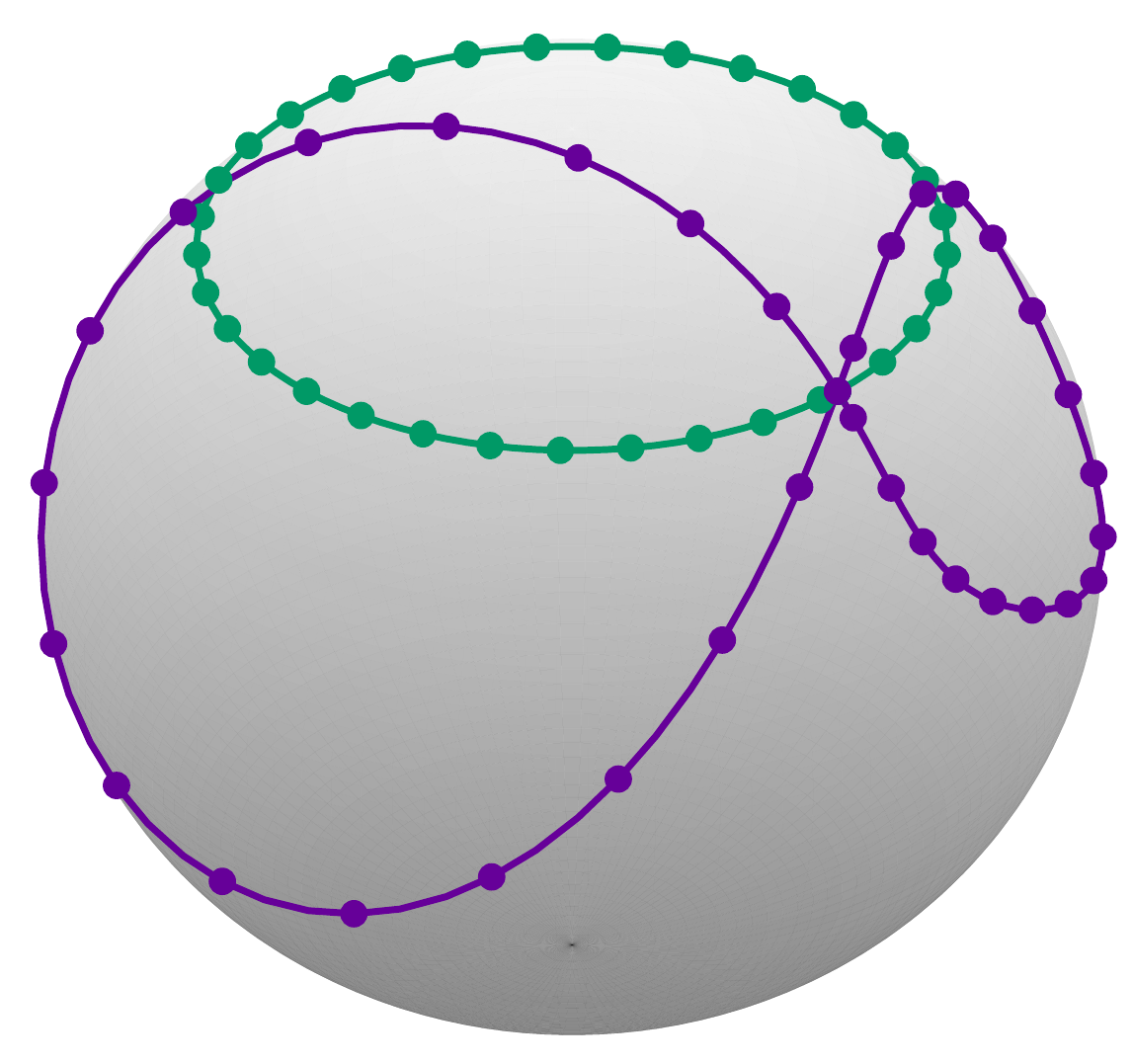}
\includegraphics[width=0.24\textwidth]{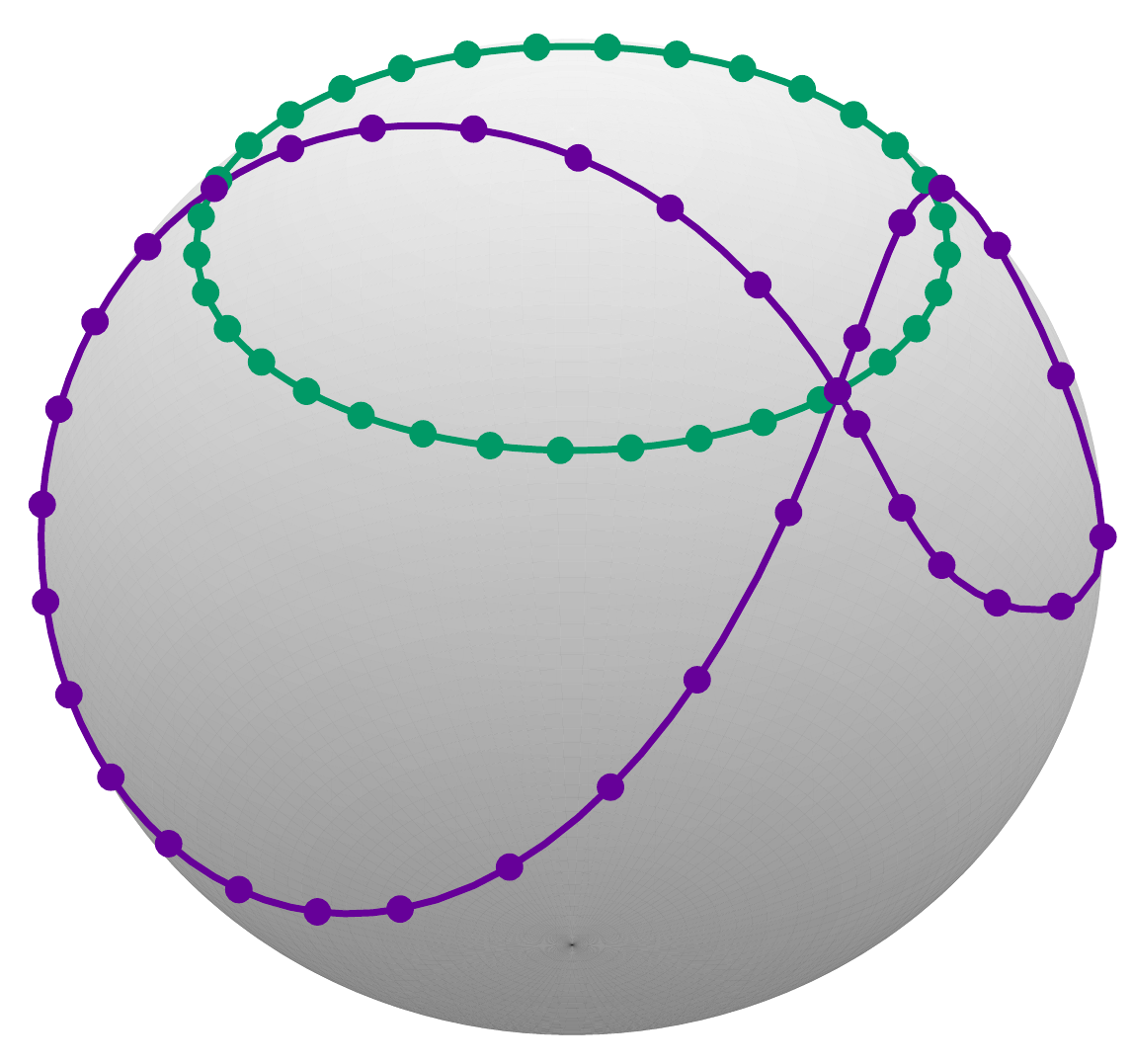}
\includegraphics[width=0.24\textwidth]{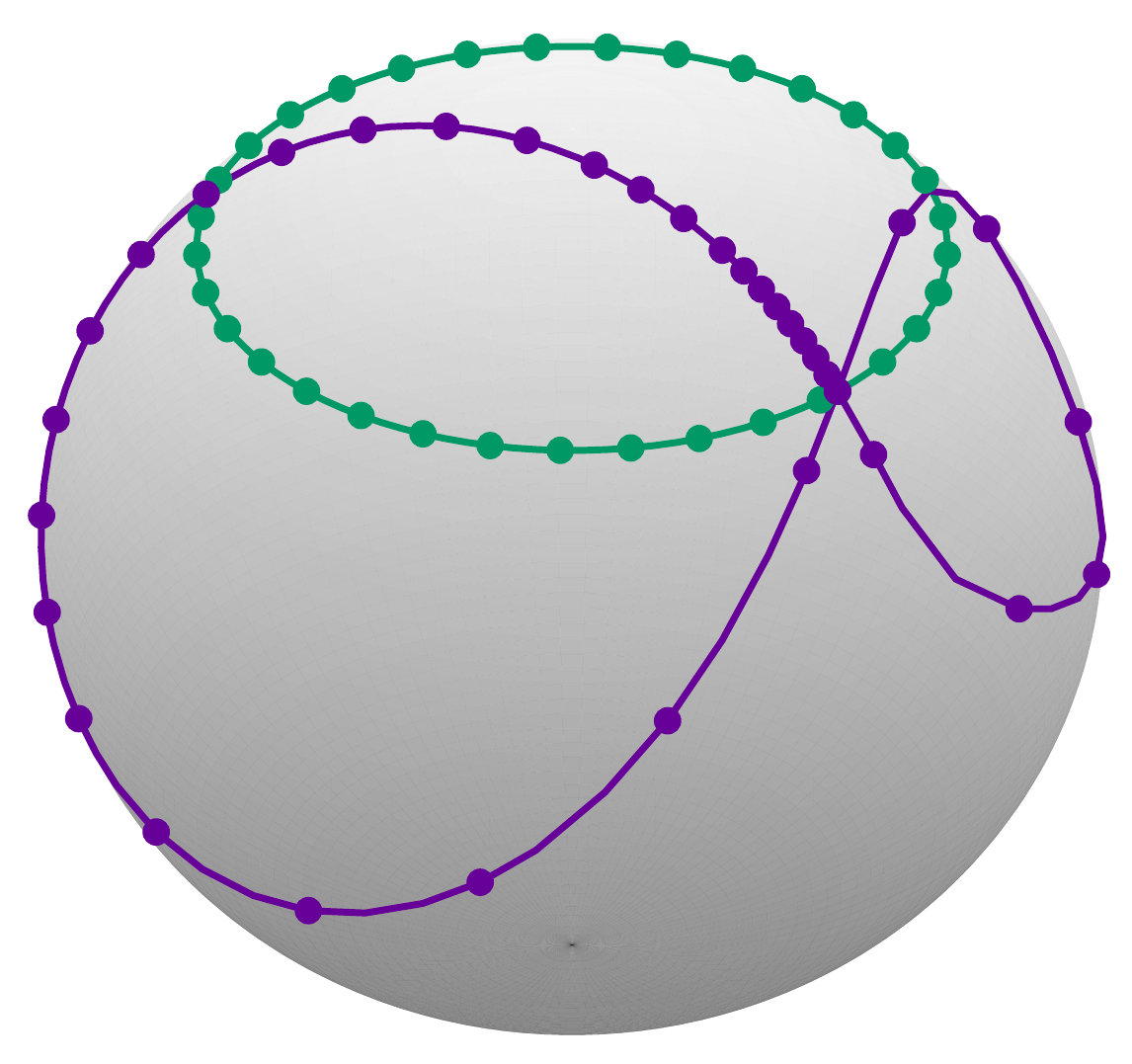}}\\
\subfloat[The interpolated curves at times $\theta = \left\{\frac{1}{4},\frac{1}{2},\frac{3}{4}\right\}$, from left to right, before reparametrization, on $\text{S}^2$ (blue line) and $\text{SO}(3)$ (yellow line).]{
\includegraphics[width=0.33\textwidth]{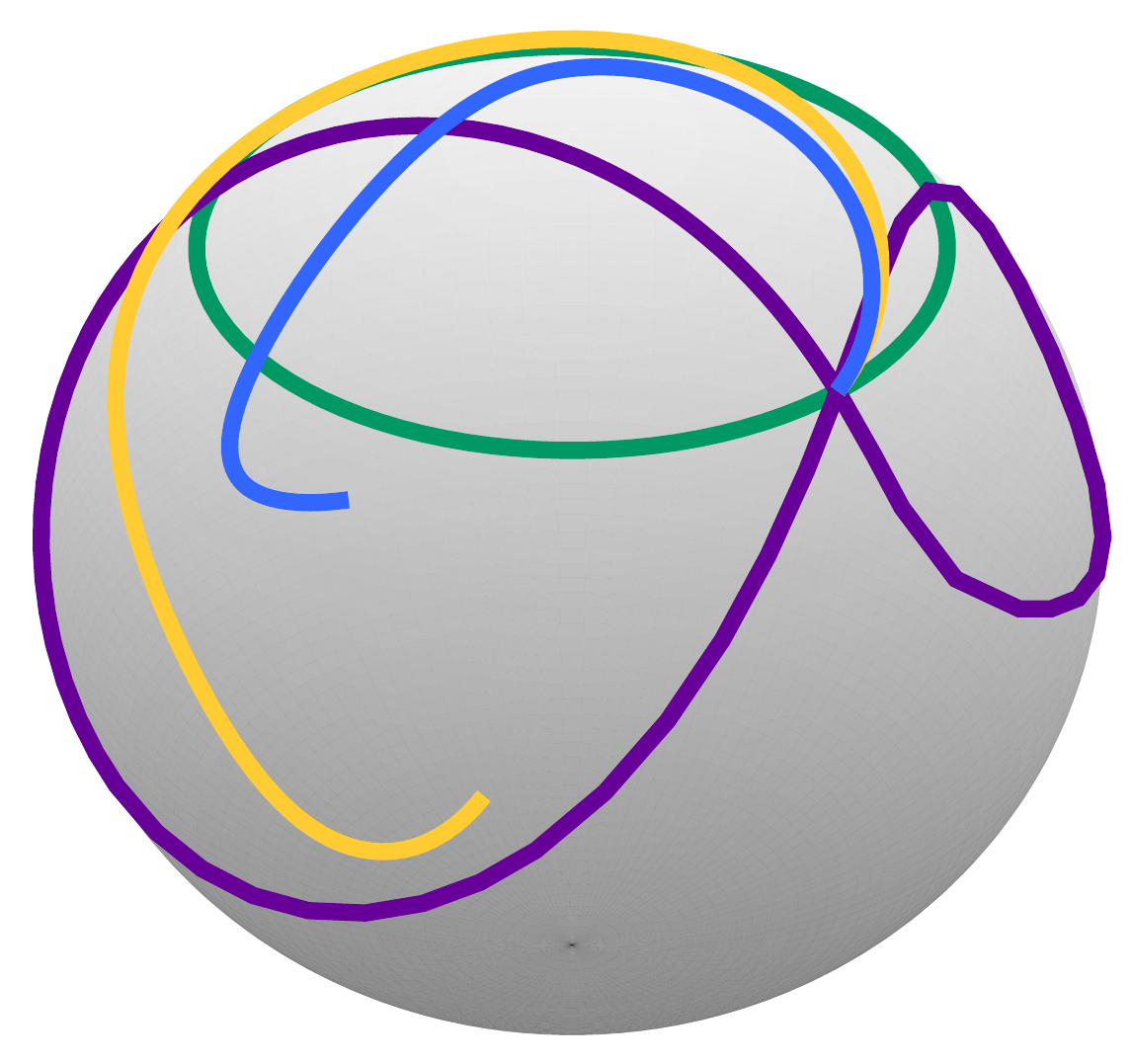}
\includegraphics[width=0.33\textwidth]{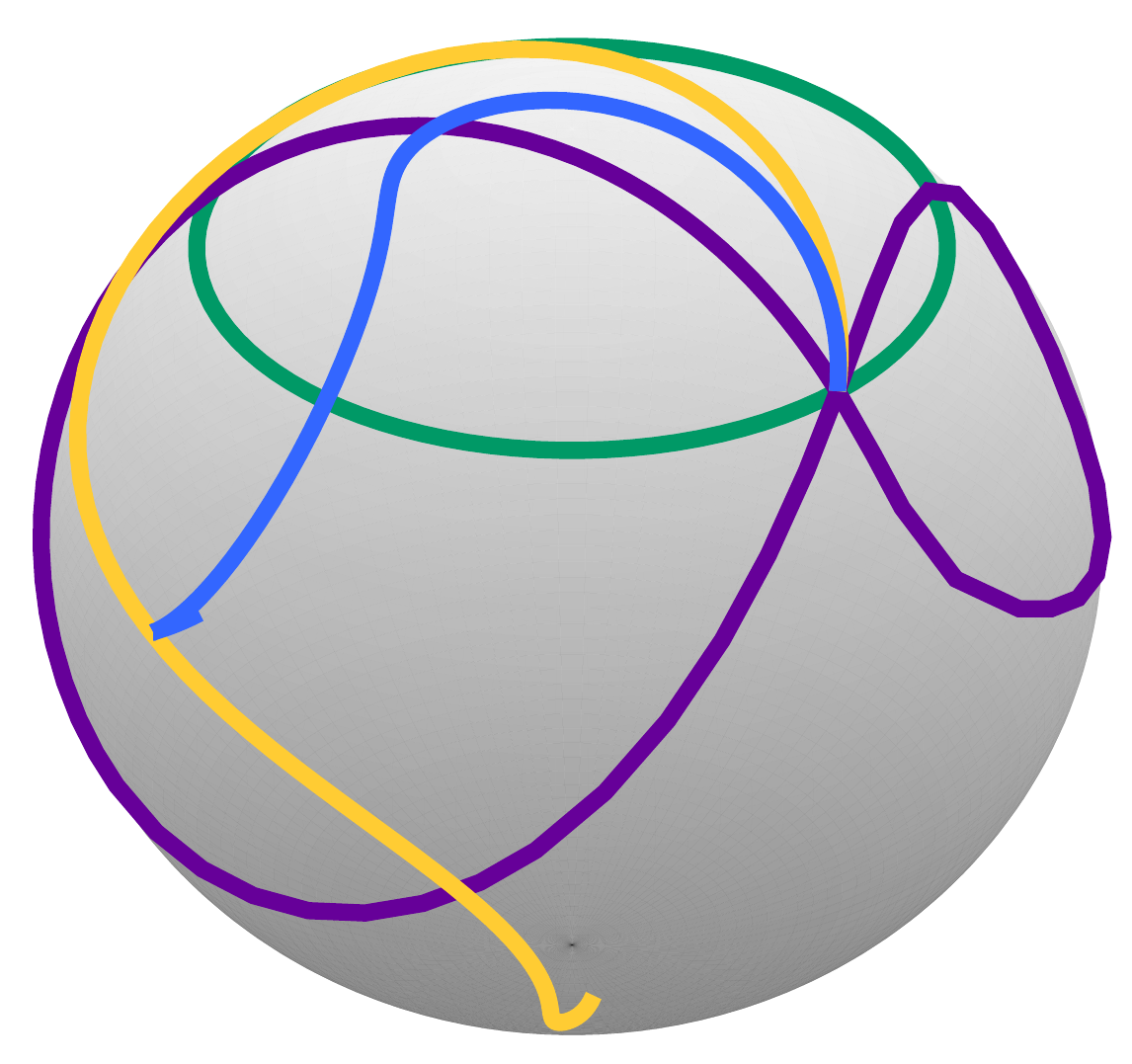}
\includegraphics[width=0.33\textwidth]{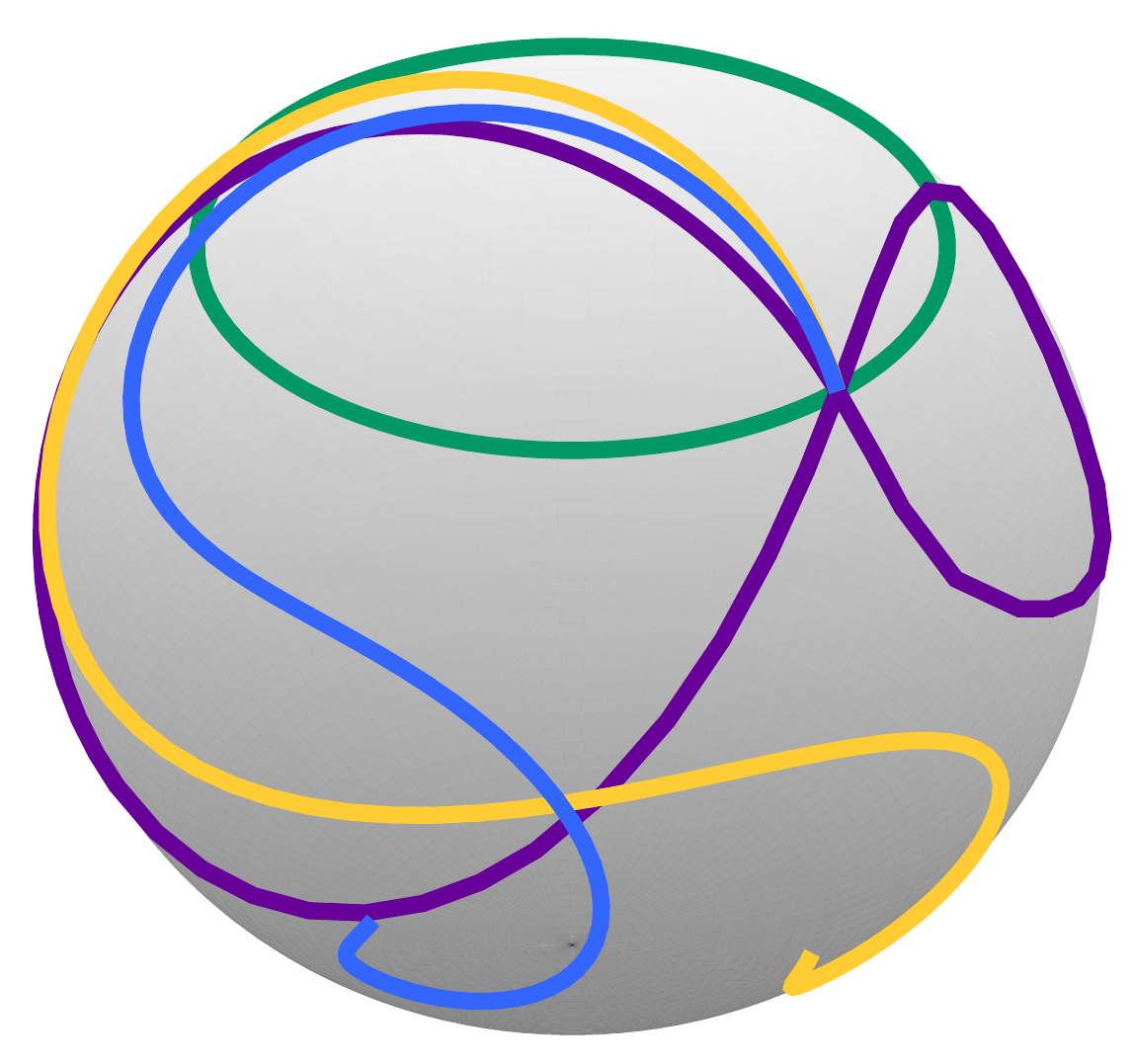}}\\
\subfloat[The interpolated curves at times $\theta = \left\{\frac{1}{4},\frac{1}{2},\frac{3}{4}\right\}$, from left to right, after reparametrization, on $\text{S}^2$ (blue line) and $\text{SO}(3)$ (yellow line).]{
\includegraphics[width=0.33\textwidth]{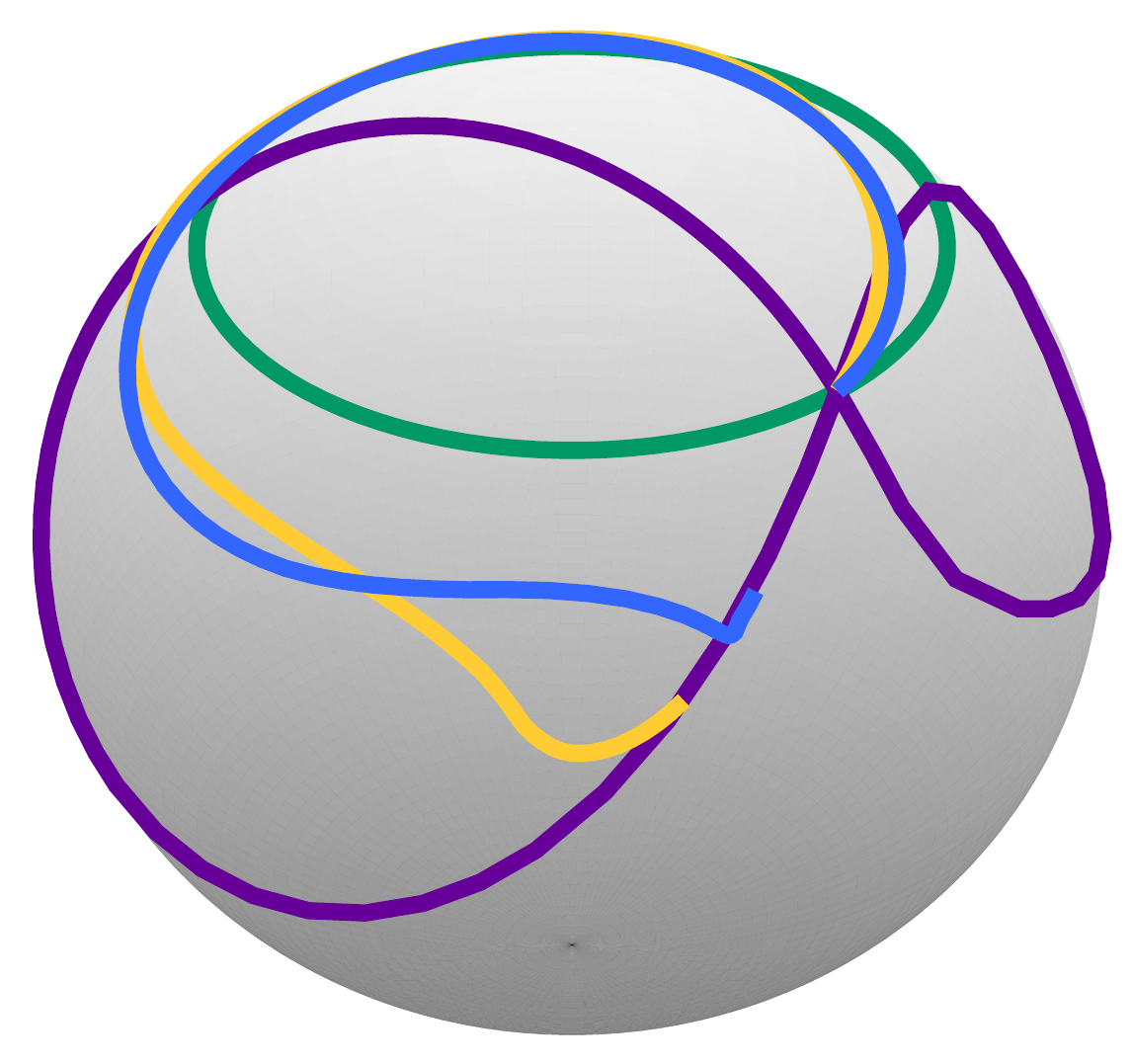}
\includegraphics[width=0.33\textwidth]{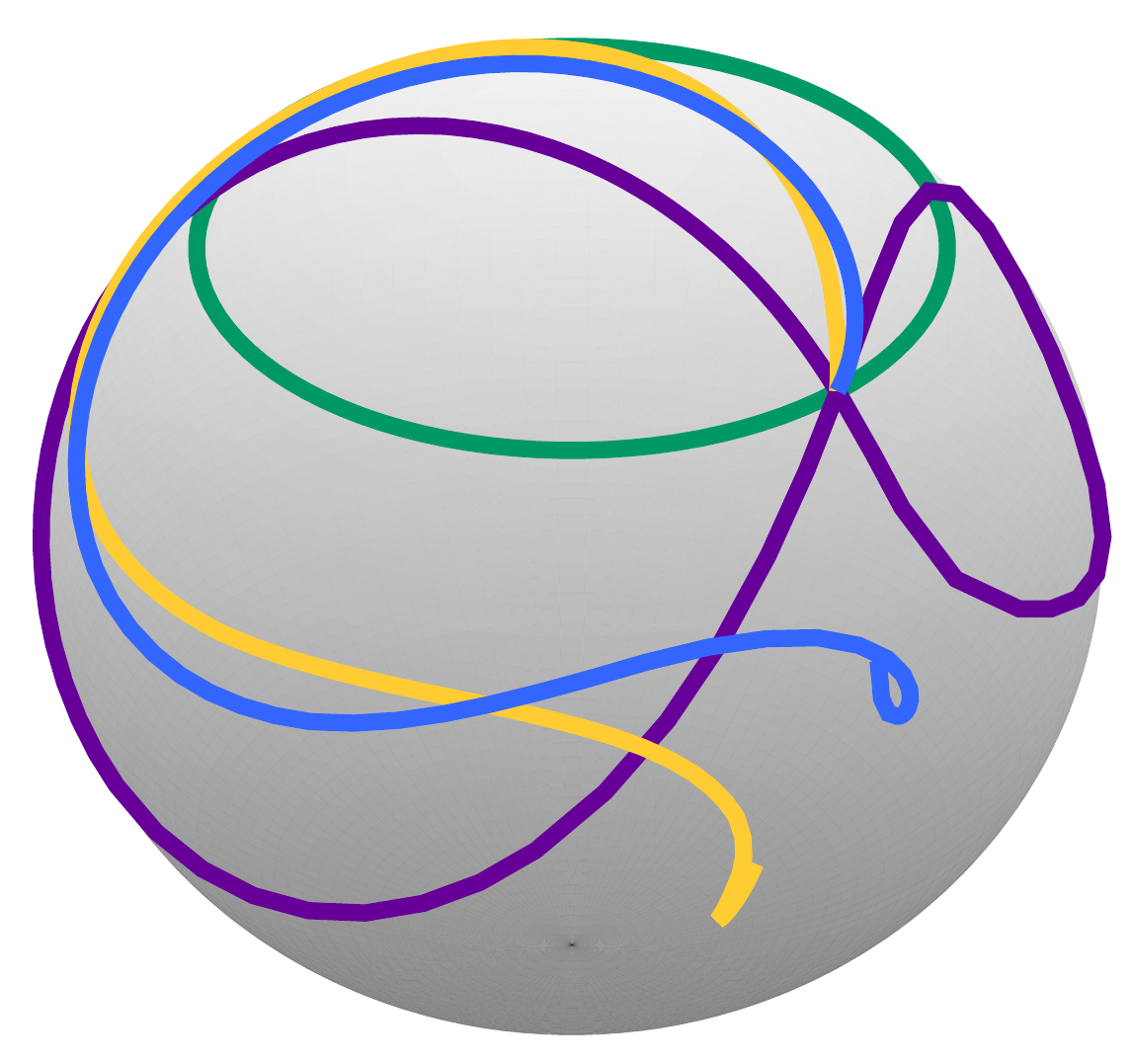}
\includegraphics[width=0.33\textwidth]{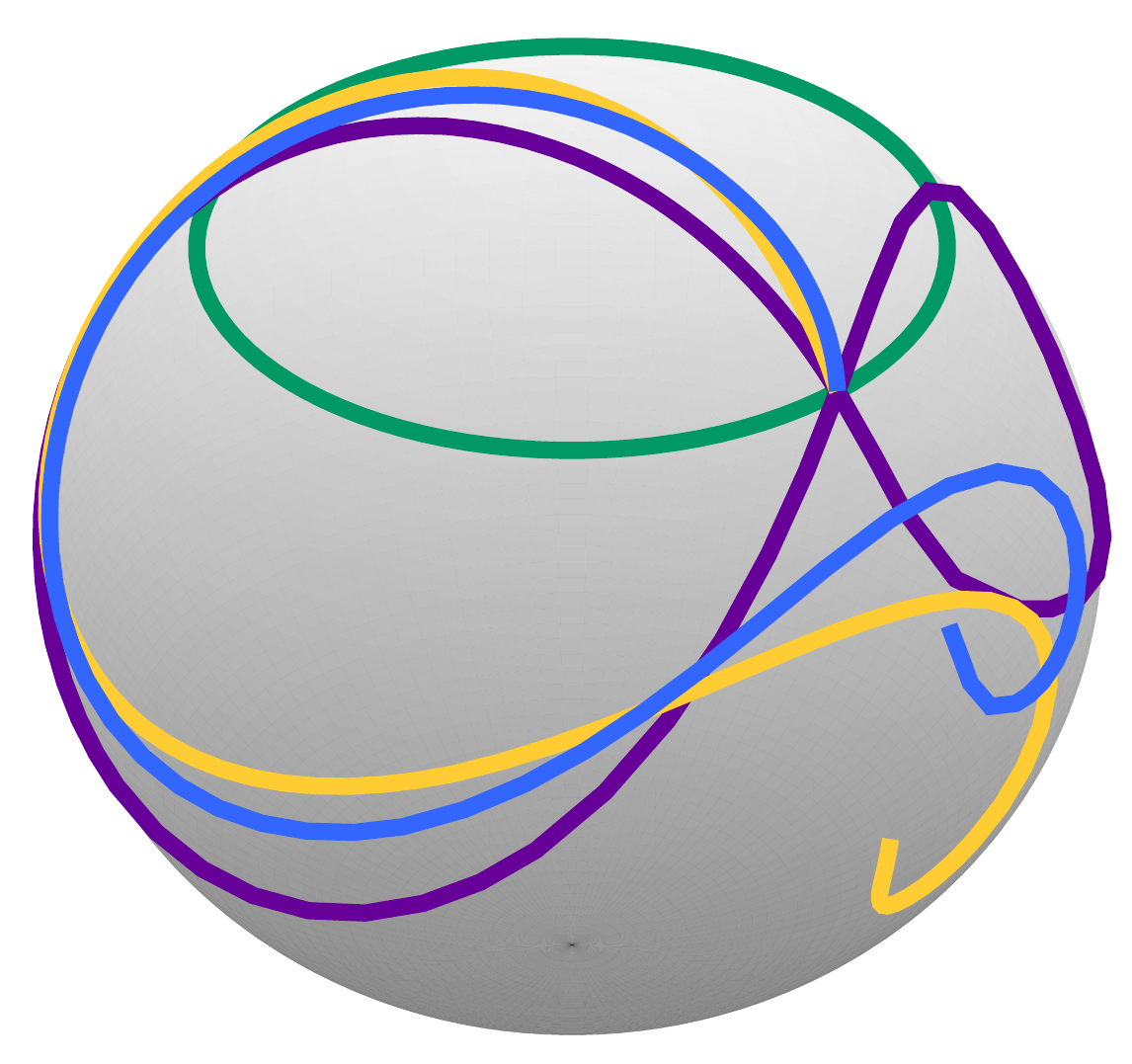}}
\caption{Interpolation between two curves on $\text{S}^2$, with and without reparametrization, obtained by the reductive SRVT. The results are compared to the corresponding SRVT interpolation between curves on $\text{SO}(3)$. The $\mathrm{SO}(3)$ curves are mapped to $\text{S}^2$ by multiplying with the vector $(0,1,1)^\text{T}/\sqrt{2}$. 
}
\label{fig:curves1}
\end{figure}
\end{center}

In Figure \ref{fig:curves1} we show instants of the computed geodesic in the shape space of curves on the sphere between two curves $\bar{c}_1$ and $\bar{c}_2$, using the reductive SRVT. We compare this to the geodesic between the curves $\bar{z}_1$ and $\bar{z}_2$ in $\mathrm{SO}(3)$ which when mapped to $\mathrm{S}^2$ gives $\bar{c}_1$ and $\bar{c}_2$. We show the results obtained before and after reparametrization. In the latter case, a dynamic programming algorithm, see \cite{sebastian03}, was used to reparametrize the curve $\bar{c}_2(t)$ such that its distance to $\bar{c}_1(t)$, measured by taking the $L^2$ norm of $\bar{q}_1(t)- \bar{q}_2(t)$ in the Lie algebra, is minimized.
The various instances of the geodesics between $\bar{c}_1(t)$ and $\bar{c}_2(t)$ are found by interpolation,
\begin{equation*}
\bar{c}_\text{int}(\bar{c}_1,\bar{c}_2,\theta) = \mathcal{R}^{-1}\left(\left(1-\theta\right)\mathcal{R}(\bar{c}_1) + \theta \,\mathcal{R}(\bar{c}_2)\right), \qquad \theta \in [0,1].
\end{equation*}

%
%

\bibliographystyle{splncs03}
\bibliography{Shape_Lie_homog}

\end{document}